\documentclass[preprint,12pt]{elsarticle}

\usepackage{amsthm}
\usepackage{amsmath}
\usepackage{tikz}
\usetikzlibrary{backgrounds,calc,cd,positioning}

\usepackage{amssymb}
\usepackage{amsmath}
\usepackage{xcolor}
\usepackage{enumitem}
\usepackage[figuresright]{rotating}

\newtheorem{theorem}{Theorem}
\newtheorem{lemma}{Lemma}
\newtheorem{defn}{Definition}
\newtheorem{rem}{Remark}
\newtheorem{ex}{Example}
\newtheorem{cor}{Corrolary}

\usepackage{amssymb}
\usepackage{amsmath}

\journal{Discrete Applied Mathematics}

\begin{document}

\begin{frontmatter}



\title{Balanced Domination in Convex Polytopes, Trees, and Grid Graphs}


\author[label1]{Bojan Nikolic\corref{cor1}}
\ead{bojan.nikolic@pmf.unibl.org}

\author[label1,label2,label3]{Marko Djukanovic}
\ead{marko.djukanovic@pmf.unibl.org}
\ead{marko.dukanovic@ung.si}

\cortext[cor1]{Corresponding author.}

\address[label1]{Faculty of Natural Sciences and Mathematics, University of Banja Luka, Mladena Stojanovića 2, Banja Luka, 78000, Bosnia and Herzegovina\fnref{label3}}

\address[label2]{University of Nova Gorica, Vipavska cesta 13, Rožna Dolina, SI-5000 Nova Gorica, Slovenia}
 
\address[label3]{Institute of Information Science (IZUM), Prešernova ulica 17, 2000 Maribor,  Slovenia\fnref{label3}}
  
\begin{abstract}
	This paper addresses two open questions posed in~\cite{xu2021balanced} regarding the balanced domination number in graphs. We show that three new classes of graphs---those of convex polytopes $A_n$, $D_n$, and $R_n''$---are $d$-balanced. Further, we provide a characterization of $d$-balancedness for rooted trees with two levels of descendants and prove that each full binary tree is $d$-balanced. Additionally, several results for caterpillar graphs are established. Moreover, we determine and prove the exact balanced domination number on grid graphs.
	Finally, we conclude the paper by providing several new  problems of interest.
\end{abstract}


\begin{keyword}
	Balanced domination number\sep $d$-balanced graphs\sep layer label-sum technique

	
	
	\MSC[2020] 05C50\sep 05C69\sep 05C70
	
	
\end{keyword}
\end{frontmatter}


\section{Introduction}\label{sec:intro}

The minimum dominating set problem and its numerous variants have attracted significant attention in both theoretical and applied graph theory~\cite{fomin2008combinatorial,megiddo1988finding,jiang2023exact,nacher2016minimum}. Formally, the domination problem is defined as follows. Given is a simple, undirected graph $G=(V,E)$, where $V$ represents the vertex set and $E \subseteq V \times V$ the edge set, the goal is to find a minimum-cardinality subset $D \subseteq V$ such that every vertex $v \in V$ either belongs to $D$ or is adjacent to at least one vertex $u \in D$, i.e., $\{u,v\} \in E$. The size of optimal set $D$ is called the \textit{domination number} of $G$. This problem has a wide range of practical applications, including bioinformatics~\cite{zhao2020minimum,nacher2016minimum}, networking~\cite{karbasi2013application}, graph mining~\cite{chalupa2018order}, and automatic text summarization~\cite{xu2016generalized}, among others. From a theoretical perspective, exact results and bounds for the domination number have been determined for many graph families, such as interval graphs~\cite{golovach2017minimal}, graph bundles~\cite{zmazek2006domination}, grid and cylindrical grid graphs~\cite{chang1992domination, gonccalves2011domination,nandi2011domination}, planar graphs~\cite{macgillivray1996domination}, chemical graphs~\cite{quadras2015domination}, Cartesian products of graphs~\cite{el1991domination}, to name a few. For comprehensive references on the fundamentals of graph domination, we refer the reader to~\cite{haynes2023domination, haynes2013fundamentals,cockayne2006domination}.

Over the past decades, numerous extensions and variations of the domination problem have been introduced and thoroughly studied, motivated by theoretical interest or notable practical applications. Among the most prominent are the strong domination number~\cite{rautenbach2000bounds}, the total domination number~\cite{lam2007total}, the connected domination number~\cite{sampathkumar1979connected}, the Roman domination number~\cite{liu2012upper,favaron2009roman}, and the $k$-domination number~\cite{kazemi2012total}, among others.

Recently, Xu et al.~\cite{xu2021balanced} introduced another related problem, the \textit{balanced domination number}, defined as follows. Let $G=(V(G),E(G))$ be a graph where each vertex can be assigned a label from $\{-1,0,1\}$. A labeling function $f$ defined on the set of vertices $V(G)$ is a \textit{balanced dominating function} (BDF) if, for every vertex $v \in V(G)$, the sum of labels across all vertices in the closed neighborhood of $v$ equals zero. Among all such functions, a BDF $f$ that maximizes the total sum $\omega_f = \sum_{v \in V(G)} f(v)$ is called the maximum balanced dominating function, and its corresponding weight $\omega_f$ is known as the \textit{balanced domination number} of $G$,  denoted by $\gamma_b(G)$. Xu et al.~\cite{xu2021balanced} established several general upper bounds on this number, some of which are tight. They introduced the concept of $d$-balanced graphs, which are those  with $\gamma_b(G)=0$. They further characterized several classes that are $d$-balanced: $r$-regular graphs, graph coronas, double stars, and complete multipartite graphs. In addition, they posed open problems for this concrete problem: 
($i$) How to characterize all $d$-balanced caterpillar graphs? and  
($ii$) How to determine the exact balanced domination number of grid graphs? 

\vspace{5pt}

The main contributions of this paper are summarized as follows:
\begin{itemize}
	\item We prove that certain graph classes of convex polytopes are $d$-balanced, specifically the families $A_n$, $D_n$, and $R_n''$. This is done based on a restated definition of the balanced domination number through the eye of linear algebra and the kernel of a linear operator. 
	
	\item We make significant progress on the open problem ($i$)   for caterpillar graphs by showing the necessary conditions on the number of leaves of a caterpillar graph that admits a labeling by a non-zero balanced domination function. In addition, we establish a characterization under which rooted trees with two levels of descendants and at least two children of the root node are $d$-balanced. On top of that, given is a proof that each full binary tree is $d$-balanced.
	
	\item We provide a complete solution to the open problem ($ii$), determining the balanced domination number for grid graphs.
\end{itemize}

\subsection{Notation and Preliminaries} \label{sec:notation}
If not stated differently, we will always deal with simple, undirected graph $G$. Let the order of graph $G$ be denoted by $|V(G)|$ (or by $|V|$ if it is clear from the context to which graph is referred to). 
For each $v \in V(G)$, let $N_G[v]=\{ u \in V(G) \colon uv \in E(G)\} \cup \{ v\} $ denotes the closed  neighborhood of vertex $v$.   For a graph $G$, let $A(G)$ 
denotes its {\it{adjacency matrix}} defined by $$A(G)=[a_{ij}]_{|V|\times |V|},\,\,		a_{ij}:=\left\{
\begin{array}{ll}
	1, & \hbox{if $v_iv_j\in E$;} \\
	0, & \hbox{otherwise.}
\end{array}
\right.$$
Thus, $a_{ij}$ indicates the presence or absence of an edge between two vertices $v_i,v_j\in V$. 

Remainder of the paper is organized as follows. Section~\ref{sec:convex_polytopes} provides the proofs that several graph classes of convex polytopes are $d$-balanced. Section~\ref{sec:trees} reveals necessary conditions for the caterpillar graphs to not be $d$-balanced and establishes a characterization of $d$-balancedness for a subclass of trees with two levels of descendants. Additionally, it has proven that full binary trees are also $d$-balanced. In Section~\ref{sec:bdn-grid-graphs} it has been shown that the grid graphs are $d$-balanced. 
Finally, in Section~\ref{sec:conclusions} conclusions and outline for future work are discussed.

\section{Balanced Domination Number on some Graphs of Convex Polytopes}\label{sec:convex_polytopes}

In this section, we give formal proofs that three classes of graphs of convex polytopes are all $d$-balanced. To this end, we first provide the formal definition of the problem and the concept of $d$-balancedness, after which we present the technique employed in establishing these results. 

\begin{defn}[\cite{xu2021balanced}] \label{defn1}
	Function $f \colon V(G) \mapsto \{-1, 0, 1\}$ is called  balanced domination function (BDF) iff for any vertex $v$  the sum of labels across all vertex in $N_{V(G)}[v]$ equals 0. The weight of a BDF $f$ is given by $\omega_f=\sum_{v \in V(G)} f(v)$. The maximum weight across all  BDF defines the Balanced domination number of $G$ ($\gamma_{bd}(G)$). 
\end{defn}

\begin{defn}[\cite{xu2021balanced}]
	A graph $G$ is called $d$-balanced, if $\gamma_{bd}(G)=0$. 
\end{defn}

Let $G$ be an arbitrary graph. Observe that the function  $f:=0$ serves as a balanced domination function on graph $G$, implying  $\gamma_{bd}(G)\geqslant 0$. Thus, to show that $G$ is a $d$-balanced graph, it suffices to prove that no balanced domination function on $G$ has a positive weight. 

Let $I(G)$ denote the identity matrix of order $|V(G)|$. We then introduce the linear operator $M(G):=A(G)+I(G)$ associated with the graph $G$. It is easy to check that the condition of a balanced dominating function, as stated in Definition \ref{defn1}, can be formulated in terms of the $M(G)$ operator. Indeed, for any BDF function $f$, and labeling $x_i=f(v_i), i\in\{1,\dotsc,|V(G)|\}$, the vector $\textbf{x}=(x_1,\dotsc,x_{|V(G)|})$ is a solution of the equation $M(G)\cdot \textbf{x}^T=\textbf{0}^T$, where $\textbf{0}$ is the $1\times |V(G)|$ zero vector. Conversely, every solution $\textbf{x}=(x_1,\dotsc,x_{|V(G)|})$ of the system of linear equations $M(G)\cdot \textbf{x}^T=\textbf{0}^T$ determines a proper BDF function $f$ defined as $f(v_i)=x_i,i\in\{1,\dotsc,|V(G)|\}$. Hence, there is
a bijective correspondence between the set of BDF functions of the graph $G$ and the kernel $ker(M(G))$ of the linear operator $M(G)$. This enables us to restate the definition of the balanced domination number of the graph $G$: $$\gamma_{bd}(G)=\max\left\{\sum_{i=1}^{|V(G)|} x_i:\textbf{x}=(x_1,\dotsc,x_{|V(G)|})\in ker(M(G))\right\}.$$
Consequently, a graph $G$ is $d$-balanced if and only if for each vector $\textbf{x}=(x_1,\dotsc,x_{|V(G)|})\in ker(M(G))$, it holds $\sum_{i=1}^{|V(G)|}x_i=0$. 

\vspace{5pt}

Convex polytope $A_n=(V(A_n),E(A_n)),n \geqslant 5$, or antiprism, was introduced in~\cite{imran2013classes}. It is a graph with the following sets of vertices
and edges
\begin{align*}
	&V(A_n) = \{u_i, v_i, w_i \mid 0\leqslant i\leqslant n-1 \},\\
	&E(A_n) = \{u_i u_{i+1}, v_iv_{i+1}, w_iw_{i+1}, u_iv_i, v_iw_i, u_{i+1}v_i, v_{i+1}w_i \mid 0\leqslant i\leqslant n-1 \}.
\end{align*}

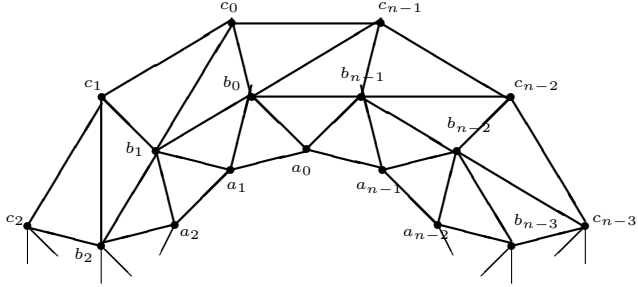
\begin{figure}	
	\def\n{7}
	\begin{center}
		\begin{tikzpicture}[scale=1,
			every node/.style={circle,draw,inner sep=0.5pt,minimum size=3.5mm},
			edge/.style={line width=0.8pt}]

			\foreach \i in {0,...,\numexpr\n-1\relax}{%
				\pgfmathsetmacro{\ang}{360*\i/\n}
				\node (u\i) at (\ang:3) {$u_{\i}$};
				\node (v\i) at (\ang:2) {$v_{\i}$};
				\node (w\i) at (\ang:1) {$w_{\i}$};
			}

			\foreach \i in {0,...,\numexpr\n-1\relax}{%
				\pgfmathtruncatemacro{\next}{mod(\i+1,\n)}
				\draw[edge] (u\i) -- (u\next);
				\draw[edge] (v\i) -- (v\next);
				\draw[edge] (w\i) -- (w\next);
				\draw[edge] (u\i) -- (v\i);
				\draw[edge] (v\i) -- (w\i);
				\draw[edge] (u\next) -- (v\i);
				\draw[edge] (v\next) -- (w\i);
			}

		\end{tikzpicture}
		\caption{The graph of convex polytope $A_7$} \label{fig:an}
	\end{center}
\end{figure}

{It is important to emphasize that the indices of the vertices of the graph $A_n$ take modulo $n$. For example, the graph of convex polytope $A_7$, is displayed in Figure~\ref{fig:an}.}
\begin{theorem}\label{thm:an_thm}
	For every $n\geqslant 5$, $A_n$ is a $d$-balanced graph.
\end{theorem}

\begin{proof}
	
	The vertices of  $A_n$ are distributed over three layers, each containing $n$ vertices. Within each layer, the vertices form a cycle of length $n$. The adjacency between the first and second layers follows the same pattern as the one between the second and third layers, while no edges connecting any vertex from  the first layer with any vertex from the third layers. Thus, the  $3n\times 3n$ format matrix of the linear operator $M(A_n)$ can be represented as a block matrix composed of the nine $n\times n$ layer-wise adjacency matrices. More precisely, the matrix of the linear operator $M(A_n)$ can be written in a block form as: 
	$$
	M(A_n) =
	\begin{bmatrix}
		C & B & O \\
		B^T & C & B \\
		O & B^T & C
	\end{bmatrix},
	$$ 
	where $C:=\begin{bmatrix}
		1 & 1 & 0 & 0 & \dotsc & 0 & 1 \\
		1 & 1 & 1 & 0 & \dotsc &0 & 0 \\
		0 & 1 & 1 & 1 & \dotsc &0 & 0 \\
		\vdots & \vdots & \vdots & \vdots &\vdots &\vdots &\vdots\\
		0 & 0 & 0 & 0 & \dotsc &1 & 1 \\
		1 & 0 & 0 & 0 & \dotsc &1 & 1
	\end{bmatrix}$ is the $n\times n$ matrix of a cycle graph, describing edges within each layer,  $B:=\begin{bmatrix}
		1 & 1 & 0 & 0 & \dotsc & 0 & 0 \\
		0 & 1 & 1 & 0 & \dotsc &0 & 0 \\
		0 & 0 & 1 & 1 & \dotsc &0 & 0 \\
		\vdots & \vdots & \vdots & \vdots &\vdots &\vdots &\vdots\\
		0 & 0 & 0 & 0 & \dotsc &1 & 1 \\
		1 & 0 & 0 & 0 & \dotsc &0 & 1
	\end{bmatrix}$ is the $n\times n$ adjacency matrix describing edges between consecutive layers, $B^T$ is its transpose matrix, while $O$ denotes the $n\times n$ zero matrix.
	
	To prove that $A_n$ is a $d$-balanced graph, it is sufficient to show that every vector $\textbf{x}=(x_1,\dotsc,x_{3n})\in ker(M(A_n))$ satisfies the condition $\sum_{i=1}^{3n} x_i=0$. 
	
	Let $\textbf{x}=(x_1,\dotsc,x_{3n})\in ker(M(A_n))$ be an arbitrary vector. We partition 
	this vector layer-wise into three subvectors  $\textbf{x}_1$, $\textbf{x}_2$ and $\textbf{x}_3$, corresponding to the layers of the graph $A_n$:
	
	subvector $\textbf{x}_1:=(x_1,\dotsc, x_n)$ corresponding to the first layer, 
	
	subvector $\textbf{x}_2:=(x_{n+1},\dotsc, x_{2n})$ corresponding to  the second layer,
	
	subvector $\textbf{x}_3:=(x_{2n+1},\dotsc, x_{3n})$ corresponding to  the third layer. 
	
	Moreover, let us define three sums of the components of  these subvectors: $s_1:=\sum_{i=1}^n x_i$, $s_2:=\sum_{i=n+1}^{2n} x_i$ and $s_3:=\sum_{i=2n+1}^{3n} x_i$. Obviously, it holds $\sum_{i=1}^{3n} x_i=s_1+s_2+s_3$.
	
	Using the previous notation, the system of the linear equations $M(A_n)\cdot \textbf{x}^T=\textbf{0}^T$ is stated as $\begin{bmatrix}
		C & B & O \\
		B^T & C & B \\
		O & B^T & C
	\end{bmatrix}\cdot \begin{bmatrix}
		\textbf{x}_1^T \\
		\textbf{x}_2^T \\
		\textbf{x}_3^T
	\end{bmatrix}=\begin{bmatrix}
		\textbf{0}^T \\
		\textbf{0}^T \\
		\textbf{0}^T
	\end{bmatrix}$, where each of the $\mathbf{0}$ on the right side is a  $n\times 1$ row vector. This leads to the system of  equations 
	\begin{align*}
		&C\cdot \textbf{x}_1^T+B\cdot \textbf{x}_2^T=\textbf{0}^T\\
		&B^T\cdot \textbf{x}_1^T+C\cdot \textbf{x}_2^T+B\cdot \textbf{x}_3^T=\textbf{0}^T\\
		&B^T\cdot \textbf{x}_2^T+C\cdot \textbf{x}_3^T=\textbf{0}^T.
	\end{align*}
	Multiply each equation in this system from the left by the $1\times n$ row vector $\textbf{1}:=(1,\dotsc,1)$. Since every column of the matrix $C$ contains exactly three ones and every column and every row of the matrix $B$ contains exactly two ones, we conclude $\textbf{1}\cdot C=3\cdot\textbf{ 1}$ and $\textbf{1}\cdot B=\textbf{1}\cdot B^T=2\cdot\textbf{ 1}$. Also, it holds $\textbf{1}\cdot \textbf{x}_i^T=s_i,i\in\{1,2,3\}$. Hence, the previous system reduces to the following system of linear equations
	\begin{align*}
		&3s_1+2s_2=0\\
		&2s_1+3s_2+2s_3=0\\
		&2s_2+3s_3=0,
	\end{align*}
	From the first equation, we have $s_1=-\frac{2}{3}s_2$ and from the third equation $s_3=-\frac{2}{3}s_2$. Substituting these values into the second equation gives $\frac{1}{3}s_2=0$, that is $s_2=0$, and consequently $s_1=s_3=0$. Therefore,  $\sum_{i=1}^{3n} x_i=s_1+s_2+s_3=0$, which concludes the proof.    	
\end{proof}
The ``layer-sum technique'' from the previous theorem extends naturally to a general framework for establishing $d$-balancedness. Let $G$ be a graph whose vertex set $V(G)$ is partitioned into $k\geqslant 2$ layers such that each block of $M(G)$ has a constant column sum $c_{ij}$. Any vector $\mathbf{x}\in\ker(M(G))$ splits layer-wise into subvectors $\mathbf{x}_1,\dotsc,\mathbf{x}_k$ with component sums $s_1,\dotsc,s_k$. Left-multiplying each block row in $M(G)\mathbf{x}^T=\mathbf{0}^T$ by the appropriate all-ones row vector yields the reduced system $\sum_{j=1}^k c_{ij}s_j=0$, $i\in\{1,\dotsc,k\}$. If every solution of this reduced system satisfies $\sum_{i=1}^k s_i=0$, then all vectors in $\ker(M(G))$ have zero component sum, proving that $G$ is $d$-balanced.

%

\vspace{5pt}

In the next two theorems, we apply the technique described above. Their proofs focus primarily on the inter-layer connection types specific to the graph classes under consideration, while omitting the routine computational details common to all cases.

The convex polytope $D_n,n\geqslant 5$, was introduced in~\cite{baca1988labelings}.  Formally, graph $D_n=(V(D_n), E(D_n))$ consists of the following sets
of vertices and edges:
\begin{align*}
	&V(D_n) = \{ u_i, v_i, w_i, q_i \mid 0\leqslant i\leqslant n-1\}, \\
	&E(D_n) = \{ u_iu_{i+1}, q_iq_{i+1},u_iv_i, v_iw_i, v_{i+1}w_i, w_iq_i  \mid 0\leqslant i\leqslant n-1\}.
\end{align*}
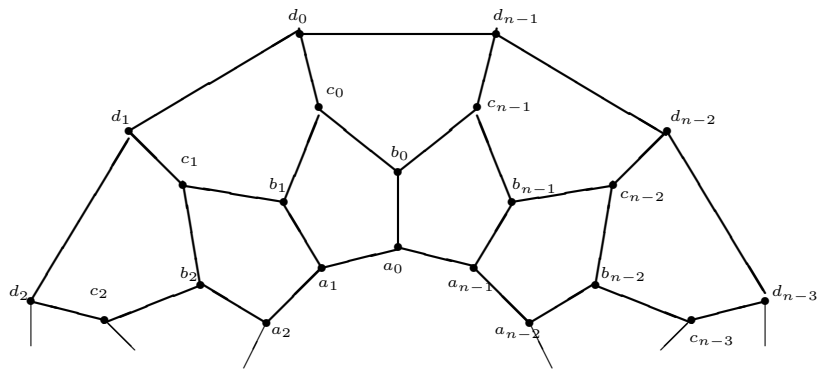
\begin{figure}
	\def\n{7}
	\begin{center}
		\begin{tikzpicture}[scale=0.7,
			every node/.style={circle,draw,inner sep=0.5pt,minimum size=3.5mm},
			edge/.style={line width=0.8pt}]
			
			\foreach \i in {0,...,\numexpr\n-1\relax}{%
				\pgfmathsetmacro{\ang}{360*\i/\n}
				\node (u\i) at (\ang:4) {$u_{\i}$};
				\node (v\i) at (\ang:3) {$v_{\i}$};
				\node (w\i) at (\ang:2) {$w_{\i}$};
				\node (q\i) at (\ang:1) {$q_{\i}$};
			}
			
			\foreach \i in {0,...,\numexpr\n-1\relax}{%
				\pgfmathtruncatemacro{\next}{mod(\i+1,\n)}
				\draw[edge] (u\i) -- (u\next);
				\draw[edge] (q\i) -- (q\next);
				\draw[edge] (u\i) -- (v\i);
				\draw[edge] (v\i) -- (w\i);
				\draw[edge] (v\next) -- (w\i);
				\draw[edge] (w\i) -- (q\i);
			}
			
		\end{tikzpicture}
		\caption{The graph of convex polytope $D_7$} \label{fig:dn}
	\end{center}
\end{figure}
We emphasize that the indices of the vertices of $D_n$
are taken modulo $n$. For example, the graph of convex polytope  $D_7$ is shown in Figure~\ref{fig:dn}.
\begin{theorem}
	For every $n\geqslant 5$, $D_n$ is a $d$-balanced graph.
\end{theorem}
\begin{proof}
	The vertices of  $D_n$ are arranged into four layers $\{u_i\},\{v_i\},\{w_i\}$ and $\{q_i\}$,  each containing $n$ vertices.
	In terms of graph connectivity, the column sums $c_{ij}$ of the matrix blocks correspond to the degrees of adjacencies within and between the layers:
		\begin{itemize}[noitemsep]
			\item 1st layer ($u_i$): $c_{11}=2+1=3$ (cycle edges $+$ self), $c_{12}=1$ ($u_i v_i$).
			\item 2nd layer ($v_i$): $c_{21}=1$ ($u_i v_i$), $c_{22}=1$ (only self), $c_{23}=2$ ($v_i w_i, v_{i+1} w_i$).
			\item 3rd layer ($w_i$): $c_{32}=2$ ($v_i w_i, v_{i+1} w_i$), $c_{33}=1$ (only self), $c_{34}=1$ ($w_i q_i$).
			\item 4th layer ($q_i$): $c_{43}=1$ ($w_i q_i$), $c_{44}=2+1=3$ (cycle edges $+$ self).
		\end{itemize}		
     Consequently,  the $4n\times 4n$ matrix of the linear operator $M(D_n)$ can be written in a block form as: 
			$$
			M(D_n) =
			\begin{bmatrix}
				C & I & O & O \\
				I & I & B & O \\
				O & B^T & I & I\\
				O & O & I & C
			\end{bmatrix},
			$$ 
			where $C,B,O$ are the same $n\times n$ matrices as in the proof of the previous theorem and $I$ is the $n\times n$ identity matrix.
		
		Let $\textbf{x}=(x_1,\dotsc,x_{4n})\in ker(M(D_n))$ be an arbitrary vector.  We partition 
			this vector layer-wise into four subvectors  
			$\textbf{x}_k:=(x_{(k-1)\cdot n+1},\dotsc, x_{k \cdot n}), k\in\{1, \ldots, 4\}$,  corresponding to the $k$-th layer of graph $D_n$, with corresponding component sum $s_k$. 		
			
			
	
	The system of linear equations $M(D_n)\cdot \textbf{x}^T=\textbf{0}^T$ and corresponding reduced system take the following form 
		\begin{align*}
			&C\cdot \textbf{x}_1^T+I\cdot \textbf{x}_2^T=\textbf{0}^T\,\,\,&3s_1+s_2=0\\
			&I\cdot \textbf{x}_1^T+I\cdot \textbf{x}_2^T+B\cdot \textbf{x}_3^T=\textbf{0}^T\,\,\,&s_1+s_2+2s_3=0\\
			&B^T\cdot \textbf{x}_2^T+I\cdot \textbf{x}_3^T+I\cdot \textbf{x}_4^T=\textbf{0}^T\,\,\,&2s_2+s_3+s_4=0\\
			&I\cdot \textbf{x}_3^T+C\cdot \textbf{x}_4^T=\textbf{0}^T,\,\,\,	&s_3+3s_4=0.
		\end{align*} 
	The unique solution of the reduced system is $(s_1,s_2,s_3,s_4)=(0,0,0,0)$. Therefore,  $\sum_{i=1}^{4n} x_i=s_1+s_2+s_3+s_4=0$, which concludes the proof.
	\end{proof}
	The graph of convex polytope $R_n''=(V(R_n''),E(R_n'')), n \geqslant 5$,  was introduced in~\cite{macdougall2006vertex}. It consists of the following sets of vertices and edges:
	{\small
		\begin{align*}
			&V(R_n'') = \{u_i, v_i, w_i, p_i, q_i, r_i \mid 0\leqslant i\leqslant n-1\}, \\
			&E(R_n'') = \{ u_i u_{i+1}, r_i r_{i+1}, u_i v_i,  v_i w_i, w_ip_i,
			p_iq_i, q_ir_i, v_{i+1} w_i, p_{i} q_{i+1}  \mid 0\leqslant i\leqslant n-1\}.
		\end{align*}
	}



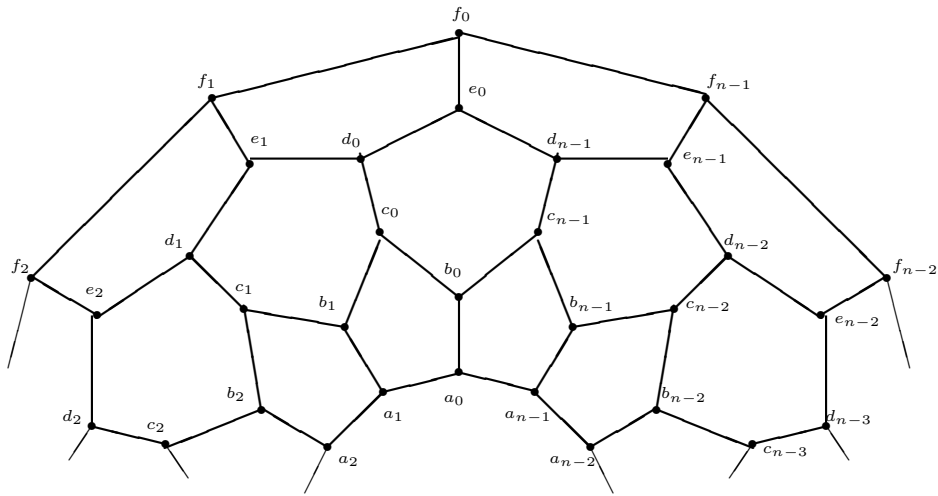
\begin{figure}[h]
	\def\n{7}
	\begin{center}	
		\begin{tikzpicture}[scale=0.65,
			every node/.style={circle,draw,inner sep=0.5pt,minimum size=3.5mm},
			edge/.style={line width=0.8pt}]

			\foreach \i in {0,...,\numexpr\n-1\relax}{%
				\pgfmathsetmacro{\ang}{360*\i/\n}
				\node (u\i) at (\ang:6) {$u_{\i}$};
				\node (v\i) at (\ang:5) {$v_{\i}$};
				\node (w\i) at (\ang:4) {$w_{\i}$};
				\node (p\i) at (\ang:3) {$p_{\i}$};
				\node (q\i) at (\ang:2) {$q_{\i}$};
				\node (r\i) at (\ang:1) {$r_{\i}$};
			}

			\foreach \i in {0,...,\numexpr\n-1\relax}{%
				\pgfmathtruncatemacro{\next}{mod(\i+1,\n)}
				\draw[edge] (u\i) -- (u\next);
				\draw[edge] (r\i) -- (r\next);
				\draw[edge] (u\i) -- (v\i);
				\draw[edge] (v\i) -- (w\i);
				\draw[edge] (v\next) -- (w\i);
				\draw[edge] (w\i) -- (p\i);
				\draw[edge] (p\i) -- (q\i);
				\draw[edge] (p\i) -- (q\next);
				\draw[edge] (q\i) -- (r\i);
			}
		\end{tikzpicture}
		\caption{The graph of convex polytope $R_7''$} \label{fig:rnsec}
	\end{center}
\end{figure}

We emphasize that the indices of the vertices of $R_n''$ are taken modulo $n$. For example, the graph of convex polytope  $R_7''$ is shown in Figure~\ref{fig:rnsec}.
\begin{theorem}
	For every $n\geqslant 5$, $R_n''$ is a $d$-balanced graph.
\end{theorem}
\begin{proof}
	 The vertices of the graph $R_n''$ are arranged into six layers $\{u_i\},\{v_i\},$ $\{w_i\}, \{p_i\}, \{q_i\}$ and $\{r_i\}$, each containing $n$ vertices. In terms of graph connectivity, the column sums $c_{ij}$ of the matrix blocks correspond to the degrees of adjacencies within and between the layers:
	 			\begin{itemize}[noitemsep]
	 					\item 1st layer ($u_i$): $c_{11}=2+1=3$ (cycle edges $+$ self), $c_{12}=1$ ($u_i v_i$).
	 					\item 2nd layer ($v_i$): $c_{21}=1$, $c_{22}=1$ (only self), $c_{23}=2$ ($v_i w_i, v_{i+1} w_i$).
	 					\item 3rd layer ($w_i$): $c_{32}=2$, $c_{33}=1$ (only self), $c_{34}=1$ ($w_i p_i$).
	 					\item 4th layer ($p_i$): $c_{43}=1$, $c_{44}=1$ (only self), $c_{45}=2$ ($p_i q_i, p_i q_{i+1}$).
	 					\item 5th layer ($q_i$): $c_{54}=2$, $c_{55}=1$ (only self), $c_{56}=1$ ($q_i r_i$).
	 					\item 6th layer ($r_i$): $c_{65}=1$, $c_{66}=2+1=3$ (cycle edges $+$ self).
	 				\end{itemize}
	 		Consequently,  the $6n\times 6n$ matrix of the linear operator $M(R_n'')$ can be written in a block form as: 
\[
	 	M(R_n'') =
	 	\begin{bmatrix}
	 		C & I & O & O & O & O \\
	 		I & I & B & O & O & O \\
	 		O & B^T & I & I & O & O\\
	 		O & O & I & I & B & O\\
	 		O & O & O & B^ T & I & I\\
	 		O & O & O & O & I & C
	 	\end{bmatrix},
\] 
	 	where $C,B,O,I$ are the same $n\times n$ matrices as in the proof of the previous theorem.
		
	Let $\textbf{x}=(x_1,\dotsc,x_{6n})\in ker(M(R_n''))$ be an arbitrary vector. We partition 
		this vector layer-wise into six subvectors  $\textbf{x}_k:=(x_{(k-1)\cdot n+1},\dotsc, x_{k \cdot n}), k\in \{1, \ldots, 6\}$ corresponding to the $k$-th layer of graph $R_n''$, with corresponding component sum $s_k$.  
	%
	%
	%
	%
	%
	%
	
	The system of linear equations $M(R_n'')\cdot \textbf{x}^T=\textbf{0}^T$ and corresponding reduced system take the following form 
			\begin{align*}
			&C\cdot \textbf{x}_1^T+I\cdot \textbf{x}_2^T=\textbf{0}^T\,\,\, &3s_1+s_2=0\\
			&I\cdot \textbf{x}_1^T+I\cdot \textbf{x}_2^T+B\cdot \textbf{x}_3^T=\textbf{0}^T\,\,\, &s_1+s_2+2s_3=0\\
			&B^T\cdot \textbf{x}_2^T+I\cdot \textbf{x}_3^T+I\cdot \textbf{x}_4^T=\textbf{0}^T\,\,\, &2s_2+s_3+s_4=0\\
			&I\cdot \textbf{x}_3^T+I\cdot \textbf{x}_4^T+B\cdot \textbf{x}_5^T =\textbf{0}^T\,\,\, &s_3+s_4+2s_5=0\\
			&B^T\cdot \textbf{x}_4^T+I\cdot \textbf{x}_5^T+I\cdot \textbf{x}_6^T=\textbf{0}^T\,\,\,&2s_4+s_5+s_6=0\\
			&I\cdot \textbf{x}_5^T+C\cdot \textbf{x}_6^T=\textbf{0}^T,\,\,\, &s_5+3s_6=0.
		\end{align*}
	The unique solution of the reduced system is $(s_1,s_2,s_3,s_4,s_5,s_6)=(0,0,0,0,$ $0,0)$. Therefore,  $\sum_{i=1}^{6n} x_i=s_1+s_2+s_3+s_4+s_5+s_6=0$, which concludes the proof.
\end{proof}
	

\section{Balanced Domination Number on some Subclasses of Trees}\label{sec:trees}

In this section we establish the characterization of $d$-balancedness of rooted trees with two levels of descendants. Moreover, some results are established for the subclass of trees consisting of caterpillar graphs.
 
\vspace{5pt}

Let $T$ be a rooted tree graph with two levels of descendants. More precisely, let $T$ have the root vertex $a_0$,
with vertices $a_1,\dotsc,a_n, n\geqslant 2$, as its children. For each $i\in\{1,\dotsc,n\}$, let $l_i$ denotes the number of children of vertex $a_i$. If $l_i\geqslant 1$, denote these vertices by $a_{i,1},\dotsc, a_{i,l_i}$. For trees of this form, the characterization of their $d$-balancedness is given in the following theorem.

\begin{theorem} \label{th:two-leveled tree}
	Let $T$ be a rooted tree graph with two levels of descendants with at least two children of the root. If $l:=\left|\left\{i\in\{1,\dotsc,n\}:l_i=2 \right\}\right|$, then $T$ is not a $d-$balanced graph if and only if $l=\frac{n-1}{2}$ and $l_i\in\{0,2\}$ holds for each $i\in\{1,\dotsc, n\}$. 
\end{theorem}

\begin{proof} Let $f$ be an arbitrary BDF function on the rooted tree graph $T$ with two levels of descendants. Define
	\begin{align*}
		&x_0:=f(a_0),\\
		&x_i:=f(a_i),i\in\{1,\dotsc,n\},\\ 
		&x_{i,j}:=f(a_{i,j}), i\in\{1,\dotsc,n\},j\in\{1,\dotsc,l_i\}.	
	\end{align*}
	Since $f$ is a balanced domination function, vertices of the tree $T$ satisfy the following system of the equations:
	{\small\begin{align*}
		&x_0+\sum_{i=1}^n x_i=0,\,\,\mbox{ condition for the root vertex } a_0, \\
		&x_0+x_i+\sum_{j=1}^{l_i} x_{i,j}=0,\,\,\mbox{ condition for the first layer vertex } a_i,i\in\{1,\dotsc,n\}),\\
		&x_i+x_{i,j}=0, \,\,\mbox{ condition for the second layer vertex } a_{i,j},j\in\{1,\dotsc,l_i\}).
	\end{align*} 
    } 
	If $l_i=0$, then $x_i=-x_0$. Moreover, from the last equation we conclude that $x_{i,j}=-x_i$ holds for each $j\in\{1,\dotsc,l_i\}$, meaning that for $l_i\geqslant 1$, it holds $x_0+x_i+(-l_i)x_i=0$, i.e. $x_0=(l_i-1)x_i$. So, for
	each $i\in\{1,\dotsc,n\}$, it holds $x_0=(l_i-1)x_i$, which means that the weight of a BDF function $f$ is given by
	\begin{align*}
		\omega_f&=x_0+\sum_{i=1}^n x_i+\sum_{i=1}^{n}\sum_{j=1}^{l_i} x_{i,j}=x_0+\sum_{i=1}^n x_i+\sum_{i=1}^n (-l_i)x_i\\
		&=x_0+\sum_{i=1}^{n}(1-l_i)x_i=x_0-nx_0=x_0(1-n).
	\end{align*} 
	
	The condition $\omega_f>0$ is equivalent to the condition $x_0\ne 0$. This, in turn, forces that for each $i\in\{1,\dotsc,n\}$, it holds $l_i\in\{0,2\}$. Indeed, if $l_i=1$, then we directly obtain $x_0=0$. On the other hand, if $l_i\geqslant 3$, then
	$|x_0|=(l_1-1)|x_i|\geqslant 2|x_i|$, which again leads to $x_0=0$, since all labels are coming from the  set $\{-1,0,1\}$.
	Furthermore, if $l_i=0$, then $x_i=-x_0$ and, if $l_i=2$, then $x_i=x_0$. In this setting, $f$ is a valid BDF function if and only if $x_0+\sum_{i=1}^n x_i=0$ holds (the balanced condition for the root of  $T$). Since $\sum_{i=1}^n x_i=(n-l)(-x_0)+lx_0=(2l-n)x_0$, the condition $x_0+\sum_{i=1}^n x_i=0$ is equivalent to equation $(1+2l-n)x_0=0$, which is satisfied if and only if $l=\frac{n-1}{2}$.
\end{proof}	

\begin{rem}
	In the preceding proof, we depart from the previously employed technique of analyzing sums across three layers of labels (the root and the two subsequent layers), even though this approach naturally suggests itself due to the evident three-layer structure of the tree under consideration. The challenge arises from the fact that the existence of a nontrivial triple $(s_1,s_2,s_3)$ of layer-sum labels satisfying the corresponding system of equations does not guarantee for a direct derivation of the “internal” condition stated in the theorem.
\end{rem}

The tree $T$ is a binary tree if every vertex has at most two children. A binary tree is a full binary tree,
if each of its vertices has exactly $0$ or $2$ children. 

\begin{theorem} \label{th:punobinarnodrvo}
	(i) If $f$ is a BDF on the full binary tree $T$ with root vertex $r$, then $f(r)=0$.
	
	(ii) Every full binary tree is $d$-balanced.
\end{theorem}

\begin{proof}
	$(i)$ Let $f$ be an arbitrary BDF on the full binary tree $T$ with root vertex $r$. Denote by $L$ the sum of all labels assigned to the leaves of  $T$ and by $I$ the sum of all labels assigned to the internal vertices (vertices with two children), excluding the root $r$. Since $f$ is a BDF on the tree $T$, for each vertex $v\in V(T)$, it holds that $\sum_{u\in N_T[v]} f(u)=0$. Observe that the label of each leaf appears in exactly two of these equations, the label of each internal vertex different from $r$ appears in exactly four of these equations and $f(r)$ appears in exactly three of these equations. Therefore, after summing up all of these equations, we get $3f(r)+2L+4I=0$, so $3f(r)=-2L-4I$. Since the right hand-side of this equation is even integer and $f(r)\in\{-1,0,1\}$, we conclude that $f(r)=0$.
	
	$(ii)$ We use part $(i)$ together with a recursive argument. Let $T$ be an arbitrary full binary tree and $f$ an arbitrary BDF function on this tree. Employing the part $(i)$, one concludes that for root $r$ of tree $T$, it holds $f(r)=0$. The key observation is that, by removing the root $r$ from the tree $T$, we obtain two full binary subtrees $T_1$ and $T_2$, such that restrictions $f_{|T_1}$ and $f_{|T_2}$ are BDFs on these subtrees. This enables us to replicate the result from part $(i)$ to conclude that the vertices that are children of $r$ must have a label zero. This ``zero-labeling'' continues recursively until all internal vertices of $T$ are exhausted and only the leaves of $T$ remain. Clearly, these leaves trivially have a label zero. Hence, $f=0$. Thus, the zero function is the only BDF on the tree $T$, so this tree must be $d$-balanced.  
\end{proof}

\begin{figure}
	\begin{center}
		\begin{tikzpicture}[
			every node/.style={circle, draw, fill=white, inner sep=2pt},
			spine/.style={thick},
			leaf/.style={thick, dashed}
			]
			
			\node (s1) at (0,0) {};
			\node (s2) [right=1.5cm of s1] {};
			\node (s3) [right=1.5cm of s2] {};
			\node (s4) [right=1.5cm of s3] {};
			\node (s5) [right=1.5cm of s4] {};
			
			\draw[spine] (s1)--(s2)--(s3)--(s4)--(s5);
			
			\node (a1) [above left=0.7cm and 0.3cm of s1] {};
			\node (a2) [below left=0.7cm and 0.3cm of s1] {};
			\draw[leaf] (s1)--(a1);
			\draw[leaf] (s1)--(a2);
			
			\node (b1) [above=0.8cm of s2] {};
			\node (b2) [above left=0.7cm and 0.3cm of s2] {};
			\node (b3) [below left=0.7cm and 0.3cm of s2] {};
			\draw[leaf] (s2)--(b1);
			\draw[leaf] (s2)--(b2);
			\draw[leaf] (s2)--(b3);
			
			
			\node (d1) [above=0.8cm of s4] {};
			\node (d2) [below=0.8cm of s4] {};
			\draw[leaf] (s4)--(d1);
			\draw[leaf] (s4)--(d2);
			
			\node (e1) [above right=0.8cm and 0.3cm of s5] {};
			\node (e2) [above=0.8cm of s5] {};
			\node (e3) [below=0.8cm of s5] {};
			\node (e4) [below right=0.8cm and 0.3cm of s5] {};
			\draw[leaf] (s5)--(e1);
			\draw[leaf] (s5)--(e2);
			\draw[leaf] (s5)--(e3);
			\draw[leaf] (s5)--(e4);
		\end{tikzpicture}
		\caption{Caterpillar graph with a spine of five vertices and 2, 3, 0, 2, and 4 leaves attached, respectively} \label{fig:caterpillar}
	\end{center}
\end{figure}
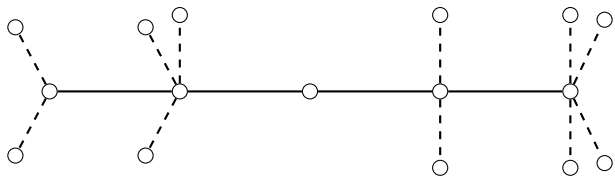

Unfortunately, it turns out that the characterization of $d$-balanced tree graphs in the general case is a considerably more challenging task. In contrast to the general case, some results can be established for a special subclass of trees consisting of caterpillar graphs. A caterpillar graph is a tree in which all vertices are within distance at most one from a central path (called the spine of the caterpillar). Equivalently, a caterpillar is a tree that becomes a simple path when all its leaves are removed. Every caterpillar graph is completely determined by the number of its vertices in the spine and by the number of leaves attached to each of these vertices. The example of caterpillar graph is shown in Fig.~\ref{fig:caterpillar}.

\vspace{5pt}

Let $C_n$ be a caterpillar graph with $n$ vertices $a_1,\dotsc,a_n$ in its spine, such that vertex $a_i$ has $l_i\geqslant 0$ leaves. For $l_i\geqslant 1$ denote these leaves by $a_{i,1},\dotsc, a_{i,l_i}$. For a given BDF function $f$ on a caterpillar graph $C_n$, let $x_i:=f(a_i), x_{i,j}:=f(a_{i,j}),i\in\{1,\dotsc,n\},j\in\{1,\dotsc,l_i\}$. Obviously, $x_{i,j}=-x_i$, for every $j\in\{1,\dotsc, l_i\}$, and, in particular, if $x_i=0$, then $x_{i,j}=0$ for all $j\in\{1,\dotsc,l_i\}$. As a consequence, for every $x_i=0$, we can replace $l_i$ with an arbitrary non negative integer without violating the balanced condition in the definition of function $f$. Specifically, if we replace all of these $l_i$'s with the value $1$, we obtain the following definition.

\begin{defn}
	Function $f \colon V(C_n) \mapsto \{-1, 0, 1\}$ is called  modified balanced domination function (MBDF) on the caterpillar graph $C_n$ iff $f$ is a BDF on the caterpillar graph $C_n$ such that $f(a_i)\in\{-1,1\}$ holds for every spine vertex $a_i\in C_n$ with  $l_i\ne 1$.	
\end{defn}

Now we can determine the number of leaves of a caterpillar graph $C_n$ that admits a non-zero MBDF. Let $L(C_n)$ be the total number of leaves of the caterpillar graph $C_n$.

\begin{theorem} \label{th:caterpillar}
	If $f$ is a non-zero MBDF on the caterpillar graph $C_n,n\geqslant 2$, then $L(C_n)\equiv (3n-2)\pmod 4$.
\end{theorem}

\begin{proof}
	Let $f$ be an arbitrary non-zero MDBF	on the caterpillar graph $C_n,n\geqslant 2$ and $x_i=f(a_i),i\in\{1,\dotsc,n\}$, represent the labels of spine vertices. Since $f$ is a non-zero MDBF, we have $x_1\ne 0,x_n\ne 0$, and  no two consecutive labels $x_i$ and $x_{i+1}$ can both be equal to zero. 
	
	Firstly, we derive explicit formulas that express $l_i$  in the terms of $x_i$ and its neighboring labels. It holds 
	\begin{align*}
		&l_1=1+x_1x_2,\\
		&l_n=1+x_{n-1}x_n,\\
		&l_i=1+x_{i-1}x_i+x_ix_{i+1}, \mbox{ for all } i\in\{2,\dotsc,n-1\}.
	\end{align*} 			
	Indeed, the balanced condition for the vertex $a_1$ yields the equation $x_1+\sum_{j=1}^{l_1} x_{i,j}+x_2=0$, which is equivalent to the equation $(1-l_1)x_1+x_2=0$. Since $x_1\in\{-1,1\}$, it holds $\frac{1}{x_1}=x_1$, and after dividing the previous equation with $x_1\ne 0$, we obtain $1-l_1+x_1x_2=0$, which immediately implies $l_1=1+x_1x_2$. Similarly, one obtains $l_n=1+x_{n-1}x_n$. For an arbitrary $i\in\{2,\dotsc,n-1\}$, if $x_i=0$, then $l_i=1$, so $l_i=1+x_{i-1}x_i+x_ix_{i+1}$ holds trivially. If $x_i\ne 0$,  the balanced condition for the vertex $a_i$ yields the equation $x_{i-1}+\sum_{j=1}^{l_i} x_{i,j}+x_{i+1}=0$, which is equivalent to the equation $x_{i-1}+(1-l_i)x_i+x_{i+1}=0$. and, since $x_i\in\{-1,1\}$, it holds $\frac{1}{x_i}=x_i$, and after dividing it with $x_i\ne 0$, we again obtain the required formula for $l_i$.
	
	Summing over all spine vertices, we obtain
	\begin{align*}
		L(C_n)&=\sum_{i=1}^{n} l_i=(1+x_1x_2)+\sum_{i=2}^{n-1}(1+x_{i-1}x_i+x_ix_{i+1})+(1+x_{n-1}x_n)\\
		&=n+2\sum_{i=1}^{n-1} x_ix_{i+1}.
	\end{align*}
	Let $S:=\sum_{i=1}^{n-1} x_ix_{i+1}$. We classify the consecutive pairs of labels $(x_i, x_{i+1})$ into three types:
	
	\begin{itemize}
		\item Type I: both labels are non-zero and have the same sign, i.e., $(1,1)$ or $(-1,-1)$. Suppose there is $p$ such pairs.
		\item Type II: both labels are non-zero and have opposite signs, i.e., $(1,-1)$ or $(-1,1)$. Suppose there is $q$ such pairs.
		\item Type III: pairs that include at least one zero. Suppose there is $r$ such pairs.
	\end{itemize}
	Obviously, it holds $p+q+r=n-1$. Additionally, since no two consecutive labels $x_i$ and $x_{i+1}$ can both be equal to zero, $r$ must be an even natural number. Consequently, $S=p-q=(p+q)-2q=n-1-r-2q$, which implies $S-(n-1)\equiv 0 \pmod 2$. Therefore, $L(C_n)=n+2((n-1)+ 2m)$, for some natural number $m$, which simplifies to
	$L(C_n)=3n-2+4m$. This concludes the proof.
\end{proof}	

\begin{rem}
	The equation $M(C_n)\cdot\textbf{x}^T=\textbf{0}^T$ is equivalent to the system of equations given by        	
	$$\begin{bmatrix}
		1-l_1 & 1 & 0 & 0 & \dotsc & 0 & 0 \\
		1 & 1-l_2 & 1 & 0 & \dotsc &0 & 0 \\
		0 & 1 & 1-l_3 & 1 & \dotsc &0 & 0 \\
		\vdots & \vdots & \vdots & \vdots &\vdots &\vdots &\vdots\\
		0 & 0 & 0 & 0 & \dotsc &1-l_{n-1} & 1 \\
		0 & 0 & 0 & 0 & \dotsc &1 & 1-l_n
	\end{bmatrix}\cdot \begin{bmatrix}
		x_1 \\
		x_2 \\
		x_3 \\
		\vdots\\
		x_{n-1}\\
		x_n
	\end{bmatrix}=\begin{bmatrix}
		0 \\
		0 \\
		0 \\
		\vdots\\
		0\\
		0
	\end{bmatrix}.$$
	In the literature is known the explicit and recursive formulas for calculating the determinant of the tridiagonal matrix on the left-hand side of this equation~\cite{qi2019some}. However, these formulas are too complicated to be used as a tool in deriving a closed-form condition for $d$-balancedness of the caterpillar graph $C_n$.      
\end{rem}

\section{Balanced Domination Number on Grid Graphs}\label{sec:bdn-grid-graphs}	

For natural numbers $m$ and $n$, the grid graph $Grid_{m\times n}$ consists of $m \cdot n$ vertices, each represented by a point $v_{ij}$ in an $m\times n$ rectangular lattice. Edges are formed between vertices that are adjacent, either horizontally within the same row or vertically within the same column. Formally, graph $Grid_{m\times n}$ can be represented as the Cartesian product $P_m\square P_n$ of two path graphs $P_m$ and $P_n$ with $m$ and $n$ vertices, respectively.

The layer-sum technique employed in the previous proofs cannot be efficiently adapted to verify the $d$-balancedness of an arbitrary grid graph. The most natural choice of layers would involve selecting either the $m$ rows or the $n$ columns of the given grid graph. However, since rows and columns of a grid graph do not have a cyclic structure, the corresponding local balanced conditions for the boundary and inner vertices take distinct forms. This makes the layer-sum approach more difficult to implement in the case of grid graphs.

For a given balanced domination function $f$ on the grid graph $Grid_{m\times n}$, let $a_{i,j}$ denote the label of the vertex $v_{i,j}$, that is, $a_{i,j}:=f(v_{i,j})$. Evidently, the knowledge of all labels in any of the four boundary lines---the first or the $m-$th row, or the first or $n-$th column---uniquely determines the function $f$.

We proceed now with showing that every grid graph is $d-$balanced. To start with, a few preliminary results will be established.
For $k\in\{1,\dotsc,m\}$, the set of labels $D_k:=\{a_{k,1},a_{k-1,2},\dotsc,a_{1,k}\}$  will be referred to as the $k-$th anti-diagonal of the given labeling.

\begin{lemma} \label{lema2} For $3\leqslant m\leqslant n$, let $f$ be a BDF on the grid graph  $Grid_{m\times n}$. If for any vertex $v\in\{v_{1,1},v_{1,n},v_{m,1},v_{m,n}\}$ holds $f\left(N[v]\right)=\{0\}$, then $f=0$.	
\end{lemma}

\begin{proof}
	Any grid graph is invariant with respect to up-down and left-to-right reflections, so it suffices to consider only the case $v=v_{1,1}$. Assume that all labels on the anti-diagonals $D_1$ and $D_2$  are equal to zero. We shall prove that each $D_k$ anti-diagonal, for $k\in\{3,\dotsc,m\}$, also consists only of zero labels. We proceed by using induction on $k$.
	
	First, consider the case $k=3$. Suppose $a_{2,2}=1$. The corresponding portion of the graph labeling is then configured as follows:
	$$\begin{tikzcd}[column sep=tiny]
		0& 0 & \boxed{-1} & \dotsc\\[-20pt]
		0 &1 &a_{2,3} & \dotsc\\[-20pt]
		\boxed{-1}& a_{3,2} & \phantom{1}& \dotsc\\[-20pt]
		\vdots & \vdots & \vdots & \vdots
	\end{tikzcd}.$$ 
	Since the function $f$ is balanced at the vertices $v_{1,3}$ and $v_{3,1}$ (labeled with the boxed -1's), we conclude $a_{3,2}\geqslant 0$ and $a_{2,3}\geqslant 0$. However, this leads to a contradiction, as it implies  $\sum_{u\in N[v_{2,2}]} f(u)\geqslant 1$, meaning $f$ is not balanced at   $v_{2,2}$. The case $a_{2,2}=-1$ is ruled out by a similar argument. Thus, $a_{2,2}=0$, and consequently, $a_{3,1}=a_{1,3}=0$.
	
	Next, for a given $k\in\{4,\dotsc,m\}$, suppose that every label above the $k-$th anti-diagonal $D_k$ is equal to zero. Furthermore, suppose that there is a non-zero label positioned "inside" this anti-diagonal. More precisely, assume that there is an index $i\in\{2,\dotsc,k-1\}$ such that $a_{k-i+1,i}\in\{-1,1\}$. Consider the case $a_{k-i+1,i}=1$. The corresponding portion of the graph labeling is then configured as follows:		
	$$\begin{tikzcd}[column sep=tiny]
		\vdots & \vdots & \vdots &\vdots\\[-20pt]
		\dotsc & 0& 0 & \boxed{-1} & \dotsc\\[-20pt]
		\dotsc & 0 &1 &a_{k-i+1,i+1} & \dotsc\\[-20pt]
		\dotsc &\boxed{-1}& a_{k-i+2,i} & \phantom{1}& \dotsc\\[-20pt]
		\vdots &\vdots & \vdots & \vdots & \vdots
	\end{tikzcd} $$
	Since the function $f$ is balanced at the vertices $v_{k-i,i+1}$ and $v_{k-i+2,i-1}$ (labeled with boxed -1's), we conclude $a_{k-i+2,i}\geqslant 0$ and $a_{k-i+1,i+1}\geqslant 0$. However, this leads to a contradiction, as it implies  $\sum_{u\in N[v_{k-i+1,i}]} f(u)\geqslant 1$, meaning that $f$ is not balanced at the vertex $v_{k-i+1,i}$. The case $a_{k-i+1,i}=-1$ is ruled out by a similar argument. Thus, for all $i\in\{2,\dotsc,k-1\}$, it holds $a_{k-i+1,i}=0$, and consequently, $a_{k,1}=a_{1,k}=0$.
	
	In this way, we have proved that there are no non-zero labels above or on the $m-$th anti-diagonal $D_m$. In particular, this means the entire first column of the grid graph is zero-labeled. Consequently, it follows $f=0$.
\end{proof}

One can easily verify that the blocks
$$ 
B_{1\times n}:=\begin{tikzcd}[column sep=5pt]
	1 & -1 &\mid & 0 & 1 &-1 & \mid & \dotsc & \mid & 0 & 1 & -1 
\end{tikzcd}
$$
$$ 
-B_{1\times n}:=\begin{tikzcd}[column sep=5pt]
	-1 & 1 &\mid & 0 & -1 &1 & \mid & \dotsc & \mid & 0 & -1 & 1 
\end{tikzcd}
$$
provide labelings for only two different non-zero BDF on the grid graph  $Grid_{1\times n}$, where $n\equiv 2\pmod 3$.
Straightforwardly, the blocks 
$$
B_{2\times n}:=\begin{tikzcd}[column sep=5pt]
	1 & \mid & 0 & -1 & \mid & 0 & 1 &\mid & 0 & -1 & \mid &\dotsc  \\[-20pt]
	-1 & \mid & 0 & 1 & \mid & 0 & -1 & \mid & 0 & 1 & \mid & \dotsc
\end{tikzcd}
$$		
$$		-B_{2\times n}:=\begin{tikzcd}[column sep=5pt]
	-1 & \mid & 0 & 1 & \mid & 0 & -1 &\mid & 0 & 1 & \mid& \dotsc  \\[-20pt]
	1 & \mid & 0 & -1 & \mid & 0 & 1 & \mid & 0 &-1 &\mid &\dotsc
\end{tikzcd}
$$
provide labelings for only two different non-zero BDF on the grid graph  $Grid_{2\times n}$, where $n\equiv 1\pmod 2$. Also, we introduce the following patterns:
\begin{center}
	$P_1:=\begin{tikzcd}[column sep=5pt]
		0 & 1 & -1 & 0\\[-20pt]
		-1& 0 & 0 & 1\\[-20pt]
		1& 0 & 0 & -1\\[-20pt]
		0 & -1 &1 & 0
	\end{tikzcd}$,
	$-P_1:=\begin{tikzcd}[column sep=5pt]
		0 & -1 & 1 & 0\\[-20pt]
		1& 0 & 0 & -1\\[-20pt]
		-1& 0 & 0 & 1\\[-20pt]
		0 & 1 &-1 & 0
	\end{tikzcd}$,
	
	$P_2:=\begin{tikzcd}[column sep=5pt]
		1 & -1 & 1 & -1\\[-20pt]
		0& -1 & 1 & 0\\[-20pt]
		0& 1 & -1 & 0\\[-20pt]
		-1 & 1 &-1 & 1
	\end{tikzcd}$,
	$-P_2:=\begin{tikzcd}[column sep=5pt]
		-1 & 1 & -1 & 1\\[-20pt]
		0& 1 & -1 & 0\\[-20pt]
		0& -1 & 1 & 0\\[-20pt]
		1 & -1 &1 & -1
	\end{tikzcd}$,
	
	$P_3:=\begin{tikzcd}[column sep=5pt]
		1 & 0 & 0 & -1\\[-20pt]
		-1& -1 & 1 & 1\\[-20pt]
		1& 1 & -1 & -1\\[-20pt]
		-1 & 0 & 0 & 1
	\end{tikzcd}$,
	$-P_3:=\begin{tikzcd}[column sep=5pt]
		-1 & 0 & 0 & 1\\[-20pt]
		1& 1 & -1 & -1\\[-20pt]
		-1& -1 & 1 & 1\\[-20pt]
		1 & 0 & 0 & -1
	\end{tikzcd}$.
\end{center}
Let $\mathbf{0}$ denotes a  zero column of format $4\times 1$. If we start at the top-left corner label of the anti-diagonal $D_1$ and, in a cascading manner, traverse through the anti-diagonals $D_k, k\in\{2,3,4\}$, excluding any labeling that would violate the balanced condition during this process, it is easy to verify that the blocks 	
\begin{align*}
	B_{4\times n}^{(t)}&=P_i\,\,|\,\,\mathbf{0}\,\,|\,\,P_t\,\,|\,\,\mathbf{0}\,\,|\,\,P_t\,\,|\,\,\dotsc\,\,|\,\,\mathbf{0}\,\,|\,\, P_t,\,\,t\in\{1,2,3\},\\
	-B_{4\times n}^{(t)}&=-P_t\,\,|\,\,\mathbf{0}\,\,|\,\,-P_t\,\,|\,\,\mathbf{0}\,\,|\,\,-P_t\,\,|\,\,\dotsc\,\,|\,\,\mathbf{0}\,\,|\,\, -P_t,\,\,t\in\{1,2,3\},
\end{align*}
enumerate all valid labelings of non-zero BDFs on the grid graph  $Grid_{4\times n}$, where $n\equiv 4\pmod 5$. In the next theorem, we show that an arbitrary nonzero labeling of the grid graph $Grid_{m\times n}$ can be constructed of using $B_{1\times n}$, $-B_{1\times n}$, $B_{2\times n}$, $-B_{2\times n}$, and $B_{4\times n}^{(t)}$, $-B_{4\times n}^{(t)},t\in\{1,2,3\}$ as the building blocks.

\begin{theorem} \label{teorema1}
	Let $f$ be a non-zero BDF on the grid graph $Grid_{m\times n}$, with $m\leqslant n$. Then, the
	labeling induced by $f$ is described by at least one of the following schemes:
	\begin{itemize}
		\item Type 1. For $m\equiv1(\mod 2)$ and $n\equiv 2(\mod 3)$, copies of blocks $B_{1\times n}$ and $-B_{1\times n}$ are alternated, with a "reset" zero row $\mathbf{0}$ of format $m\times 1$ inserted between each of them. Given that $m$ is odd, the zero-row cannot appear as the last row in this arrangement:
		\begin{center}
			\begin{tabular}{c|c}
				\mbox{Type } 1.1 & \mbox{Type } 1.2 \\
				\hline
				$B_{1\times n}$ & $-B_{1\times n}$\\
				$\mathbf{0}$ & 	$\mathbf{0}$\\
				$-B_{1\times n}$& $B_{1\times n}$\\
				$\mathbf{0}$& $\mathbf{0}$\\
				$B_{1\times n}$ & $-B_{1\times n}$\\
				$\vdots$& 	$\vdots$\\
				$\mathbf{0}$& $\mathbf{0}$\\
				$(-1)^{\frac{m-1}{2}}B_{1\times n}$& 	$(-1)^{\frac{m+1}{2}}B_{1\times n}$	    	
			\end{tabular}	
		\end{center}
		\item Type 2. For $m\equiv2(\mod 3)$ and $n\equiv 1(\mod 2)$, copies of block $B_{2\times n}$ (respectively block $-B_{2\times n}$) are merged, with a "reset" zero row $\mathbf{0}$ of format $m\times 1$ inserted between each of them. Given that $3\nmid m$, the zero-row cannot appear as the last row in this arrangement:
		\begin{center}
			\begin{tabular}{c|c}
				\mbox{Type } 2.1 & \mbox{Type } 2.2 \\
				\hline
				$B_{2\times n}$ & $-B_{2\times n}$\\
				$\mathbf{0}$ & 	$\mathbf{0}$\\
				$B_{2\times n}$& $-B_{2\times n}$\\
				$\mathbf{0}$& $\mathbf{0}$\\
				$B_{2\times n}$ & $-B_{2\times n}$\\
				$\vdots$& 	$\vdots$\\
				$\mathbf{0}$& $\mathbf{0}$\\
				$B_{2\times n}$& 	$-B_{2\times n}$	    	
			\end{tabular}	
		\end{center}
		\item Type 3. For $m\equiv4(\mod 5)$ and $n\equiv 4(\mod 5)$, copies of block $B_{4\times n}^{(i)},i\in\{1,2,3\}$, (respectively block $-B_{4\times n}^{(i)},i\in\{1,2,3\}$) introduced in the proof of the Lemma \ref{lema2} are merged, with a "reset" zero row $\mathbf{0}$ of format $m\times 1$ inserted between each of them. Given that $5\nmid m$, the zero-row cannot appear as the last row in this arrangement:
		\begin{center}
			\begin{tabular}{c|c|c}
				\mbox{Type } 3.1 & \mbox{Type } 3.2 \\
				\hline
				$B_{4\times n}^{(t)}$ & $-B_{4\times n}^{(t)}$\\
				$\mathbf{0}$ & 	$\mathbf{0}$\\
				$B_{4\times n}^{(t)}$& $-B_{4\times n}^{(t)}$\\
				$\mathbf{0}$& $\mathbf{0}$ & $t\in\{1,2,3\}$\\
				$B_{4\times n}^{(t)}$ & $-B_{4\times n}^{(t)}$\\
				$\vdots$& 	$\vdots$\\
				$\mathbf{0}$& $\mathbf{0}$\\
				$B_{4\times n}^{(t)}$& 	$-B_{4\times n}^{(t)}$	    	
			\end{tabular}	
		\end{center}
	\end{itemize}	
\end{theorem}

\begin{proof}
	In the proof of the previous lemma, we showed that a zero BDF on the grid graph is fully determined by zero labels along its anti-diagonals $D_1$ and $D_2$. This was established inductively, via a zigzag traversal of the graph along its anti-diagonals, starting from  $D_3$ and concluding with $D_m$.  The same anti-diagonal argument will be applied here to prove  that any non-zero BDF on the grid graph $Grid_{m\times n}$ is uniquely determined by its labels on anti-diagonals $D_1,D_2$ and $D_3$.
	
	Set $a_{i,j}:=f(v_{i,j})$. Clearly, the labels on the anti-diagonals $D_1$ and $D_2$ satisfy the condition $a_{2,1}+a_{1,2}=-a_{1,1}$. Employing the local balance condition, for $k\in\{3,\dotsc,m\}$ we obtain  the following relations for the labels along the anti-diagonal $D_k$:
	{\footnotesize
		\begin{align} 
		&a_{k,1}+a_{k-1,2}=-(a_{k-1,1}+a_{k-2,1}),\nonumber\\
		&a_{1,k}+a_{2,k-1}=-(a_{1,k-1}+a_{1,k-2}),\tag{$BC_k$} \label{eq:rel}\\
		&a_{k-i,i+1}+a_{k-i-1,i+2}=-(a_{k-i-1,i+1}+a_{k-i-1,i}+a_{k-i-2,i+1}),k\geqslant 4,i\in\{1,\dotsc,k-3\}\nonumber.
	\end{align}
    }
	For $k\in\{2,\dotsc,m\}$ and $l\in\{0,\dotsc,\left\lfloor \frac{k}{2}\right\rfloor-1\}$, define $s_l(k):=a_{k-l,l+1}+a_{l+1,k-l}$. Clearly, $s_0(2)=a_{2,1}+a_{1,2}=-a_{1,1}$. We shall prove that every term $s_l(k)$, with $k\in\{3,\dotsc,8\}$, $l\in\{0,\dotsc,\lfloor \frac{k}{2}\rfloor-1\}$, can be expressed inductively using the values $s_{l'}(k')$, for some $k'\in\{2,\dotsc,k-1\}$ and $l'\in\{0,\dotsc,\lfloor \frac{k'}{2}\rfloor-1\}$. Indeed, by summing up the appropriate equations from (\ref{eq:rel}), we obtain the following recurrence relations:
	{\scriptsize
	\begin{align*}
		D_3: &s_0(3)=-2a_{2,2}-2a_{1,1}-s_0(2),\\
		D_4: &s_1(4)=-a_{2,2}-s_0(2),\,\, s_0(4)+2s_1(4)=-a_{2,2}-s_0(3)-2s_0(2),\\
		D_5: &s_1(5)=-2a_{3,3}-2a_{2,2}-s_1(4)-s_0(3), \\
		&s_0(5)+2s_1(5)=-2a_{3,3}-2a_{2,2}-s_0(4)-s_1(4)-2s_0(3),\\
		D_6: &s_2(6)=-a_{3,3}-s_1(4), s_1(6)+2s_2(6)=-a_{3,3}-s_1(5)-s_0(4)-2s_1(4),\\
		&s_0(6)+2s_1(6)+2s_2(6)=-a_{3,3}-s_0(5)-s_1(5)-2s_0(4)-2s_1(4),\\
		D_7: &s_2(7)=-2a_{4,4}-2a_{3,3}-s_2(6)-s_1(5), \\
		&s_1(7)+2s_2(7)=-2a_{4,4}-2a_{3,3}-s_1(6)-s_2(6)-s_0(5)-2s_1(5),\\
		&s_0(7)+2s_1(7)+2s_2(7)=-2a_{4,4}-2a_{3,3}-s_0(6)-s_1(6)-s_2(6)-2s_0(5)-2s_1(5),\\
		D_8: &s_3(8)=-a_{4,4}-s_2(6), s_2(8)+2s_2(8)=-a_{4,4}-s_2(7)-s_1(6)-2s_2(6),\\
		&s_1(8)+2s_2(8)+2s_3(8)=-a_{4,4}-s_1(7)-s_2(7)-s_0(6)-2s_1(6)-2s_2(6),\\
		&s_0(8)+2s_1(8)+2s_2(8)+2s_3(8)=-a_{4,4}-s_0(7)-s_1(7)-s_2(7)-2s_0(6)-2s_1(6)-2s_2(6).
	\end{align*}
    }
	Using the initial condition $s_0(2)=a_{2,1}+a_{1,2}=-a_{1,1}$, the explicit forms of the above sums are:
	{\footnotesize
	\begin{align*}
		&a_{2,1}+a_{1,2}=s_0(2)=-a_{1,1},\,\,\,a_{3,1}+a_{1,3}=s_0(3)=-a_{1,1}-2a_{2,2},\\
		&a_{4,1}+a_{1,4}=s_0(4)=a_{1,1}+3a_{2,2},\,\,\, a_{3,2}+a_{2,3}=s_1(4)=a_{1,1}-a_{2,2},\\
		&a_{5,1}+a_{1,5}=s_0(5)=-2a_{2,2}+2a_{3,3},\,\,\, a_{4,2}+a_{2,4}=s_1(5)=a_{2,2}-2a_{3,3},\\
		&a_{6,1}+a_{1,6}=s_0(6)=3a_{2,2}-5a_{3,3},\,\,\, a_{5,2}+a_{2,5}=s_1(6)=-a_{1,1}-4a_{2,2}+3a_{3,3},\\ &a_{4,3}+a_{3,4}=s_2(6)=-a_{1,1}+a_{2,2}-a_{3,3}, a_{7,1}+a_{1,7}=s_0(7)=-a_{2,2}-2a_{4,4},\tag{$SC$} \label{eq:sums}\\
		&a_{6,2}+a_{2,6}=s_1(7)=3a_{2,2}+2a_{4,4}\,\,\, a_{5,3}+a_{3,5}=s_2(7)=a_{1,1}-a_{2,2}-2a_{4,4},\\
		&a_{8,1}+a_{1,8}=s_0(8)=-a_{1,1}+5a_{2,2}+7a_{4,4},\,\,\, a_{7,2}+a_{2,7}=s_1(8)=a_{1,1}-2a_{2,2}-5a_{4,4},\\
		&a_{6,3}+a_{3,6}=s_2(8)=2a_{2,2}+3a_{4,4},\,\,\, a_{5,4}+a_{4,5}=s_3(8)=a_{1,1}+a_{2,2}.
	\end{align*}
    }
	Since all labels $a_{i,j}$ come from the set $\{-1,0,1\}$, the condition $a_{6,1}+a_{1,6}=3a_{2,2}-5a_{3,3}$
	implies $a_{3,3}=a_{2,2}$, while the condition $a_{8,1}+a_{1,8}=s_0(8)=-a_{1,1}+5a_{2,2}+7a_{4,4}$ implies $a_{4,4}=-a_{2,2}$. Additionally, the condition $a_{4,1}+a_{1,4}=s_0(4)=a_{1,1}+3a_{2,2}$ imposes the restrictions $a_{1,1}+a_{2,2}\in\{-1,0,1\}$ and $a_{1,1}=0\implies a_{2,2}=0$. Hence, all sums in (\ref{eq:sums}) can be expressed solely in terms of labels $a_{1,1}$ and $a_{2,2}$.
	
	The local balance conditions together with the set of conditions~(\ref{eq:sums}) impose quite restrictive constraints on the existence of non-zero BDF labelings on the $D_1-D_8$ anti-diagonals of the given grid graph. Firstly, if we fix $a_{1,1}$, $a_{2,2}$, and one of the labels $a_{3,1}$ or $a_{1,3}$, there are exactly five configurations for the first three anti-diagonals $D_1$, $D_2$, and $D_3$ (up to the reflection symmetry $f\mapsto -f$): 
	
	$C_1:\begin{tikzcd}[column sep=-7pt]
		1 & -1 & 0\\[-20pt]
		0& 0 &  & \\[-20pt]
		-1&  & &[-20pt]
	\end{tikzcd}$, 
	$C_2:\begin{tikzcd}[column sep=-7pt]
		1 & 0 & -1\\[-20pt]
		-1& 0 &  & \\[-20pt]
		0&  & &[-20pt]
	\end{tikzcd}$,
	$C_3:\begin{tikzcd}[column sep=-7pt]
		1 & 0 & 0\\[-20pt]
		-1& -1 &  & \\[-20pt]
		1&  & &[-20pt]
	\end{tikzcd}$,
	$C_4:\begin{tikzcd}[column sep=-7pt]
		1 & -1 & 1\\[-20pt]
		0& -1 &  & \\[-20pt]
		0&  & &[-20pt]
	\end{tikzcd}$,
	$C_5:\begin{tikzcd}[column sep=-7pt]
		0 & 1 & -1\\[-20pt]
		-1& 0 &  & \\[-20pt]
		1&  & &[-20pt]
	\end{tikzcd}$.

	Furthermore, the imposed conditions uniquely determine all labels populating the anti-diagonals $D_4$–$D_8$. The key observation that, concluding with the $D_6$ anti-diagonal for configuration $C_1$, the $D_7$ anti-diagonal for configuration $C_2$, and the $D_8$ anti-diagonal for configurations $C_3$, $C_4$, and $C_5$, the defining pattern starts to repeat, effectively resetting the labeling to that of the $D_1$, $D_2$ and $D_3$ anti-diagonals. This periodic appearance then continues, thereby ensuring the uniqueness of the configuration across the entire grid graph $Grid_{m\times n}$. We shall verify  this fact in detail for the configurations $C_5$, whereas a similar consideration may be applied to the remaining configurations. 
	
	For configuration $C_5$, we set $a_{1,1}=a_{2,2}=0$, $a_{2,1}=a_{1,3}=-1$ and $a_{1,2}=a_{3,1}=1$. Then $a_{4,1}+a_{1,4}=0=a_{3,2}+a_{2,3}$ and $a_{3,3}=0$. Taking $a_{2,3}=1$ leads to the configuration 
	{\small$\begin{tikzcd}[column sep=-2pt]
		0 & 1 & -1 & \boxed{-1}\\[-20pt]
		-1& 0 &  1& 0 \\[-20pt]
		1& -1 & 0&\\[-20pt]
		1& & & 
	\end{tikzcd}$}, which is not admissible, since vertex $v_{1,4}$ (labeled with the boxed $-1$) violates the local balance condition. On the other hand, taking $a_{2,3}=-1$ leads to the configuration
    {\small 
	$\begin{tikzcd}[column sep=-2pt]
		0 & 1 & -1 & 1\\[-20pt]
		-1& 0 &  \boxed{-1}& 0 \\[-20pt]
		1& 1 & 0&\\[-20pt]
		-1& & & 
	\end{tikzcd}$}, which is also not admissible, since vertex $v_{2,3}$ (labeled with the boxed $-1$) violates the local balance condition. So, $a_{2,3}=0$ remains the only possibility, which implies $a_{3,2}=a_{4,1}=a_{1,4}=0$. Next, all label on the anti-diagonal $D_5$ are forced, leading to the configuration
    {\small 
	$\begin{tikzcd}[column sep=-2pt]
		0 & 1 & -1 & 0 & 0\\[-20pt]
		-1& 0 & 0 & 1 \\[-20pt]
		1& 0 & 0&\\[-20pt]
		0& -1& &\\[-20pt]
		0 & & & 
	\end{tikzcd}$}. The labels on the anti-diagonal $D_6$ satisfy  $a_{6,1}+a_{1,6}=a_{5,2}+a_{2,5}=a_{4,3}+a_{3,4}=0$. Moreover, we have $a_{4,4}=0$. One of the labels $a_{5,2}$ and $a_{3,4}$ must be equal to $-1$. Taking $a_{2,5}=-1$ gives the configuration 
    {\footnotesize
	$\begin{tikzcd}[column sep=-2pt]
		0 & 1 & -1 & 0 & 0 & \boxed{1}\\[-20pt]
		-1& 0 & 0 & 1 & -1 & 1\\[-20pt]
		1& 0 & 0& 0 & -1\\[-20pt]
		0& -1& 0& 0\\[-20pt]
		0 &1 & &\\[-20pt]
		-1 
	\end{tikzcd}$}, which is not admissible, as vertex $v_{1,6}$ (labeled with the boxed $1$) violates the local balance condition. Hence, $a_{2,5}=0$, and consequently $a_{1,6}=a_{6,1}=a_{5,2}=0$, $a_{4,3}=1$ and $a_{3,4}=-1$. Using the "central" label $a_{4,4}$, we iteratively determine each label on the anti-diagonal $D_7$, leading to the configuration
    {\footnotesize 
	$\begin{tikzcd}[column sep=-2pt]
		0 & 1 & -1 & 0 & 0 & 0 & 1\\[-20pt]
		-1& 0 & 0 & 1 & 0 & -1\\[-20pt]
		1& 0 & 0& -1 & 0\\[-20pt]
		0& -1& 1& 0\\[-20pt]
		0 &0 & 0&\\[-20pt]
		0 & 1\\[-20pt]
		-1
	\end{tikzcd}$}. The labels on the anti-diagonal $D_8$ satisfy $a_{8,1}+a_{1,8}=a_{7,2}+a_{2,7}=a_{6,3}+a_{3,6}=a_{5,4}+a_{4,5}=0$. One of the labels $a_{1,8}$ and $a_{2,7}$ must be equal to $-1$. Taking $a_{1,8}=0$ leads to the configuration
    {\footnotesize
	$\begin{tikzcd}[column sep=-2pt]
		0 & 1 & -1 & 0 & 0 & 0 & 1 & 0\\[-20pt]
		-1& 0 & 0 & 1 & 0 & -1 & \boxed{-1}\\[-20pt]
		1& 0 & 0& -1 & 0\\[-20pt]
		0& -1& 1& 0\\[-20pt]
		0 &0 & 0&\\[-20pt]
		0 & 1\\[-20pt]
		-1
	\end{tikzcd}$}, which is not admissible, since vertex $v_{2,7}$ (labeled with the boxed $-1$) violates the local balance condition. Hence, $a_{1,8}=-1$, implying $a_{8,1}=1$, and iteratively forcing further $a_{7,2}=a_{5,4}=a_{4,5}=a_{2,7}=0$, $a_{6,3}=-1$ and $a_{3,6}=1$. Therefore, the only possible outcome that can be acquired using the "triangle" labels from $C_5$ is
    {\footnotesize
	$\begin{tikzcd}[column sep=-2pt]
		0 & 1 & -1 & 0 & 0 & 0 & 1 & -1\\[-20pt]
		-1& 0 & 0 & 1 & 0 & -1 & 0\\[-20pt]
		1& 0 & 0& -1 & 0 & 1\\[-20pt]
		0& -1& 1& 0 & 0\\[-20pt]
		0 &0 & 0& 0\\[-20pt]
		0 & 1 & -1\\[-20pt]
		-1 & 0\\[-20pt]
		1
	\end{tikzcd}$}. In this scheme, the upper $4\times 4$ submatrix precisely correspond to the pattern $P_1$. Adjacent to this pattern, there appear a zero column on its right and a zero row below it. Beyond these, the labeling resumes the same triangular arrangement as in configuration $C_5$, forming two symmetric “triangles” — one extending to the right of the zero column and the other below the zero row. This ensures that pattern $P_1$ periodically repeats in both the rows and columns. Since labeling in the last row and column cannot terminate with zero row and zero column, we conclude that $m\equiv4\pmod 5$ and $n\equiv 4\pmod 5$. This proves that configuration $C_1$ yields non-zero BDF labeling of type 3.
	
	Using this approach for remaining configurations, we conclude that configuration $C_1$ yields a non-zero BDF labeling of type 1 with $m\equiv1\pmod 2$ and $n\equiv 2\pmod 3$, configuration $C_2$ yields a non-zero BDF labeling of type 2 with $m\equiv2\pmod 3$ and $n\equiv 1\pmod 2$ and configurations $C_3$ and $C_4$ yields a non-zero BDF labelings of type 3 with $m\equiv4\pmod 5$ and $n\equiv 4\pmod 5$.
\end{proof}

\begin{ex} There is no non-zero BDF on the Grid graph $Grid_{6\times n}$, because $6\not\equiv1(mod\ 2)$, $6\not\equiv2 (mod\ 3)$ and $6\not\equiv4 (mod\  5)$. On the other hand, either of the following two arrangements
	\begin{center}
		\begin{tabular}{c|c}
			$B_{2\times 17}$ & $B_{1\times 17}$\\
			$\mathbf{0}$ & 	$\mathbf{0}$\\
			$B_{2\times 17}$& $-B_{1\times 17}$\\
			& $\mathbf{0}$\\
			& $B_{1\times 17}$ 	    	
		\end{tabular}	
	\end{center}
	with zero row $\mathbf{0}$ of format $1\times 17$, determines a non-zero BDF on a grid graph $Grid_{5\times 17}$.	
\end{ex}

\begin{cor}
	Every grid graph $Grid_{m\times n}$, with $3\leqslant m\leqslant n$, is $d-$balanced.
\end{cor}	

\begin{proof}
	Let $f$ be an arbitrary non-zero BDF on the grid graph $Grid_{m\times n}$, with $3\leqslant m\leqslant n$. According to Theorem  \ref{teorema1}, the labeling induced by $f$ is determined by at least one of the  three given schemes. Each  scheme utilizes $B_{1\times n}$, $B_{2\times n}$, $-B_{1\times n}$,
	$-B_{2\times n}$, $B_{4\times n}^{(i)}$, $-B_{4\times n}^{(i)}$, $i\in\{1,2,3\}$, as building blocks, arranged  such that any two consecutive blocks are separated by a zero-row. This structure ensures that each block is isolated from the rest of the graph labeling.  Since each of these blocks is $d-$ balanced, the sum of weights within each block equals $0$, implying that $\omega_f=0$. This holds for an arbitrary BDF on the given grid graph, leading to the conclusion $\gamma_{bd}(Grid_{m\times n})=0$. 
\end{proof}

\section{Conclusions and Future Work} \label{sec:conclusions}

This work investigated the \emph{balanced domination number} (BDN), addressing two open questions originally posed in~\cite{xu2021balanced} which originally introduced the concept of this number. In particular, we proved that three classes of graphs corresponding to convex polytopes, namely $A_n$, $D_n$, and $R_n''$, are $d$-balanced, meaning their BDN equals zero.  The proving technique relies on a \emph{layer label-sum technique}, which proved sufficient for demonstrating the $d$-balancedness of the aforementioned graph classes. Furthermore, we established that grid graphs are $d$-balanced, thereby resolving the first open question. Concerning the second open question---we derived necessary conditions for the existence of a non-zero modified $d$-balanced domination function on \emph{caterpillar graphs}, as a special subclass of trees. In addition, we provided a characterization under which  rooted trees with two levels of descendants and at least two children of the root are $d$-balanced. As a notable consequence, we proved that every full binary tree is $d$-balanced.  

\medskip

Future research should aim to identify additional classes of graphs for which the layer label-sum technique offers a sufficient approach to establishing $d$-balancedness. This technique appears to be particularly effective for graphs with symmetrical edge structures, though a formal analysis of this observation remains undertaken.  

Moreover, $d$-balanced conditions for trees and caterpillar graphs merit further investigation in a more general framework. Based on the results presented in Theorems~\ref{th:two-leveled tree} and~\ref{th:caterpillar}, the following open questions arise:  

\begin{enumerate}
	\item \textit{Layer-based verification:}  
	Let $G$ be a graph that can be partitioned into layers such that each layer exhibits symmetrical edge structures.  
	Does the layer label-sum technique described in this paper suffice to determine whether $G$ is $d$-balanced?  
	In other words, can the $d$-balancedness of $G$ always be established solely through layer-sum considerations, without reference to other intrinsic graph properties not captured by this approach?  
	
	\item \textit{Binary trees:}  
	What necessary and sufficient conditions must be satisfied for a binary tree to be $d$-balanced?  
	
	\item \textit{Caterpillar substructure analysis:}  
	Let $f$ be a non-zero modified balanced domination function on a caterpillar graph $C_n$.  
	The labeling induced by $f$ partitions $C_n$ into subgraph blocks $C_{n_1}, \dotsc, C_{n_k}$, with $n = n_1 + \dotsc + n_k$, where each restriction $f_{|C_{n_i}}$ is a proper modified balanced domination function on $C_{n_i}$.  
	Is it possible to describe all distinct ways in which such a partition of $C_n$ can be made ``unbalanced'', in the sense that the arrangement of the blocks $C_{n_1}, \dotsc, C_{n_k}$ yields a non-zero value of $\omega_f$?  
	If so, what kind of combinatorial argument should be imposed to further constrain the condition $L(C_n) \equiv (3n - 2) \pmod{4}$?  
\end{enumerate}

\section*{Acknowledgments}

Bojan Nikolic partially funded by the Ministry of Science and Technology Development and Higher Education of the Republic of Srpska (project no.\ 1259114). This publication is co-funded by the European Union’s Horizon
Europe research and innovation program under the Marie Sklodowska-Curie COFUND
Postdoctoral Programme grant agreement No.101081355– SMASH and by the Republic of Slovenia and the European Union from the European Regional Development
Fund. Co-funded by European Union. Views and opinions expressed are those of the
authors only and do not necessarily reflect those of the European Union or European
Research Executive Agency (REA). Neither the European Union nor the REA can be
held responsible for them.

\bibliographystyle{elsarticle-num}
\bibliography{bib}

\begin{thebibliography}{10}
\expandafter\ifx\csname url\endcsname\relax
  \def\url#1{\texttt{#1}}\fi
\expandafter\ifx\csname urlprefix\endcsname\relax\def\urlprefix{URL }\fi
\expandafter\ifx\csname href\endcsname\relax
  \def\href#1#2{#2} \def\path#1{#1}\fi

\bibitem{xu2021balanced}
B.~Xu, W.~Sun, S.~Li, C.~Li, On the balanced domination of graphs, Czechoslovak
  Mathematical Journal 71~(4) (2021) 933--946.

\bibitem{fomin2008combinatorial}
F.~V. Fomin, F.~Grandoni, A.~V. Pyatkin, A.~A. Stepanov, Combinatorial bounds
  via measure and conquer: Bounding minimal dominating sets and applications,
  ACM Transactions on Algorithms (TALG) 5~(1) (2008) 1--17.

\bibitem{megiddo1988finding}
N.~Megiddo, U.~Vishkin, On finding a minimum dominating set in a tournament,
  Theoretical Computer Science 61~(2-3) (1988) 307--316.

\bibitem{jiang2023exact}
H.~Jiang, Z.~Zheng, An exact algorithm for the minimum dominating set problem.,
  in: IJCAI, 2023, pp. 5604--5612.

\bibitem{nacher2016minimum}
J.~C. Nacher, T.~Akutsu, Minimum dominating set-based methods for analyzing
  biological networks, Methods 102 (2016) 57--63.

\bibitem{zhao2020minimum}
D.~Zhao, G.~Xiao, Z.~Wang, L.~Wang, L.~Xu, Minimum dominating set of multiplex
  networks: Definition, application, and identification, IEEE Transactions on
  Systems, Man, and Cybernetics: Systems 51~(12) (2020) 7823--7837.

\bibitem{karbasi2013application}
A.~H. Karbasi, R.~E. Atani, Application of dominating sets in wireless sensor
  networks, Int. J. Secur. Its Appl 7 (2013) 185--202.

\bibitem{chalupa2018order}
D.~Chalupa, An order-based algorithm for minimum dominating set with
  application in graph mining, Information Sciences 426 (2018) 101--116.

\bibitem{xu2016generalized}
Y.-Z. Xu, H.-J. Zhou, Generalized minimum dominating set and application in
  automatic text summarization, in: Journal of Physics: Conference Series, Vol.
  699, IOP Publishing, 2016, p. 012014.

\bibitem{golovach2017minimal}
P.~A. Golovach, P.~Heggernes, M.~M. Kant{\'e}, D.~Kratsch, Y.~Villanger,
  Minimal dominating sets in interval graphs and trees, Discrete Applied
  Mathematics 216 (2017) 162--170.

\bibitem{zmazek2006domination}
B.~Zmazek, J.~Žerovnik, On domination numbers of graph bundles, Journal of
  Applied Mathematics and Computing 22~(1) (2006) 39--48.

\bibitem{chang1992domination}
T.~Y. Chang, Domination numbers of grid graphs, University of South Florida,
  1992.

\bibitem{gonccalves2011domination}
D.~Gon{\c{c}}alves, A.~Pinlou, M.~Rao, S.~Thomass{\'e}, The domination number
  of grids, SIAM Journal on Discrete Mathematics 25~(3) (2011) 1443--1453.

\bibitem{nandi2011domination}
M.~Nandi, S.~Parui, A.~Adhikari, The domination numbers of cylindrical grid
  graphs, Applied Mathematics and Computation 217~(10) (2011) 4879--4889.

\bibitem{macgillivray1996domination}
G.~MacGillivray, K.~Seyffarth, Domination numbers of planar graphs, Journal of
  Graph Theory 22~(3) (1996) 213--229.

\bibitem{quadras2015domination}
J.~Quadras, A.~Sajiya Merlin~Mahizl, I.~Rajasingh, R.~Sundara~Rajan, Domination
  in certain chemical graphs, Journal of Mathematical Chemistry 53~(1) (2015)
  207--219.

\bibitem{el1991domination}
M.~El-Zahar, C.~Pareek, Domination number of products of graphs, Ars Combin 31
  (1991) 223--227.

\bibitem{haynes2023domination}
T.~W. Haynes, S.~T. Hedetniemi, M.~A. Henning, Domination in graphs: Core
  concepts, Springer, 2023.

\bibitem{haynes2013fundamentals}
T.~W. Haynes, S.~Hedetniemi, P.~Slater, Fundamentals of domination in graphs,
  CRC press, 2013.

\bibitem{cockayne2006domination}
E.~Cockayne, Domination of undirected graphs—a survey, in: Theory and
  Applications of Graphs: Proceedings, Michigan May 11--15, 1976, Springer,
  2006, pp. 141--147.

\bibitem{rautenbach2000bounds}
D.~Rautenbach, Bounds on the strong domination number, Discrete Mathematics
  215~(1-3) (2000) 201--212.

\bibitem{lam2007total}
P.~C.~B. Lam, B.~Wei, et~al., On the total domination number of graphs,
  Utilitas Mathematica 72 (2007) 223.

\bibitem{sampathkumar1979connected}
E.~Sampathkumar, H.~Walikar, The connected domination number of a graph, J.
  Math. Phys 13~(6) (1979).

\bibitem{liu2012upper}
C.-H. Liu, G.~J. Chang, Upper bounds on roman domination numbers of graphs,
  Discrete Mathematics 312~(7) (2012) 1386--1391.

\bibitem{favaron2009roman}
O.~Favaron, H.~Karami, R.~Khoeilar, S.~M. Sheikholeslami, On the roman
  domination number of a graph, Discrete Mathematics 309~(10) (2009)
  3447--3451.

\bibitem{kazemi2012total}
A.~Kazemi, On the total k-domination number of graphs, Discussiones
  Mathematicae Graph Theory 32~(3) (2012) 419--426.

\bibitem{imran2013classes}
M.~Imran, A.~Ahmad, A.~Semani{\v{c}}ov{\'a}-fe{\v{n}}ov{\v{c}}{\'\i}kov{\'a},
  et~al., On classes of regular graphs with constant metric dimension, Acta
  Mathematica Scientia 33~(1) (2013) 187--206.

\bibitem{baca1988labelings}
M.~Baca, Labelings of 2 classes of convex polytopes, Utilitas Mathematica 34
  (1988) 24--31.

\bibitem{macdougall2006vertex}
J.~A. MacDougall, M.~Miller, M.~Ba{\v{c}}a, Vertex-magic total labeling of
  generalized petersen graphs and convex polytopes, Journal of Combinatorial
  Mathematics and Combinatorial Computing 58 (2006) 89--99.

\bibitem{qi2019some}
F.~Qi, W.~Wang, D.~Lim, B.-N. Guo, Some formulas for determinants of
  tridiagonal matrices in terms of finite generalized continued fractions, in:
  Nonlinear Analysis: Problems, Applications and Computational Methods,
  Proceedings of the 6th International Congress of the Moroccan Society of
  Applied Mathematics (SM2A 2019) organized by Sultan Moulay Slimane
  University, Facult{\'e} des Sciences et Techniques, BP, Vol. 523, 2019.

\end{thebibliography}

\end{document}